\documentclass[conference,english,usletter,final]{IEEEtran}

\usepackage{mathrsfs}
\usepackage{graphicx}
\usepackage{amsmath}
\usepackage{amssymb}
\usepackage{enumerate}
\usepackage{cmll}
\usepackage{lmodern}

\newtheorem{lemma}{Lemma}

\newtheorem{remark}{Remark}

\DeclareMathOperator{\PG}{PG}

\DeclareMathOperator{\rank}{rk}

\newcommand{\mat}[1]{\mathbf{#1}}
\newcommand{\gauss}[3]{\genfrac{[}{]}{0pt}{}{#1}{#2}_{#3}}

 \newcommand{\smax}{\mathrm{A}}
 \newcommand{\imat}{\mathbf{I}}
 \newcommand{\trace}{\mathrm{Tr}}
 
\newcommand{\rdist}{\mathrm{d}_{\mathrm{r}}}
\newcommand{\sdist}{\mathrm{d}_{\mathrm{s}}}
\newcommand{\F}{\mathbb{F}}

\title{A New Approach to the Main Problem of Subspace Coding}

\author{\IEEEauthorblockN{Liu Haiteng} \IEEEauthorblockA{Department of
    Information Science\\
    and Electronics Engineering\\
    Zhejiang University, 38 Zheda Road\\
    310027 Hangzhou, China\\
    Email: liuhaiteng@zju.edu.cn} \and
  
\thanks{This work
    was supported by the National Natural Science Foundation of China
    (Grant No.\ 60872063)}

  \IEEEauthorblockN{Thomas Honold} \IEEEauthorblockA{Department of
    Information Science\\
    and Electronics Engineering\\
    Zhejiang University, 38 Zheda Road\\
    310027 Hangzhou, China\\
    Email: honold@zju.edu.cn}}
  
\begin{document}

\maketitle

\begin{abstract}
  
  Subspace codes form the appropriate mathematical setting for
  investigating the Koetter-Kschischang model of fault-tolerant
  network coding. 
  The Main Problem of
  Subspace Coding asks for the determination of a subspace code of
  maximum size (proportional to the transmission rate) if the
  remaining parameters are kept fixed. We describe a new approach to
  finding good subspace codes, which surpasses the known size limit of
  lifted MRD codes and is capable of yielding an alternative
  construction of the currently best known binary subspace code of
  packet length $7$, constant dimension $3$ and minimum subspace
  distance $4$.  
\end{abstract}


\section{Introduction}\label{sec:intro}

Let $q>1$ be a prime power. A $q$-ary (constant-dimension) subspace
code with parameters $(v,M,d;k)$ is a set
$\mathcal{C}=\{U_1,\dots,U_M\}$ of $M$ distinct $k$-dimensional
subspaces $U_i$ of the ``ambient'' vector space $\F_q^v$ (over $\F_q$)
with minimum subspace distance
$\sdist(\mathcal{C})=\min\bigl\{\sdist(U_i,U_j);1\leq i<j\leq
M\bigr\}=d$. Here the \emph{subspace distance} of $U$ and $V$ is
defined in general as $\sdist(U,V)=\dim(U+V)-\dim(U\cap V)$ and
in the special case $\dim(U)=\dim(V)=k$ reduces to $\sdist(U,V)=2k-2\dim(U\cap
V)$. In particular, $d=2\delta$ is
always even and $t=k-\delta+1$ is the smallest integer such that every
$t$-dimensional subspace of $\F_q^v$ is contained in at most one
member of $\mathcal{C}$.

The \emph{Main Problem of (constant-dimension) Subspace Coding} can be
described as follows:

\begin{quote}
  Given a prime power $q>1$ and positive integers $v,k,\delta$ with
  $2\leq\delta\leq k\leq v/2$, determine the largest cardinality $M$
  of a $q$-ary $(v,M,2\delta;k)$ subspace code. This cardinality will
  be denoted by $\smax_q(v,2\delta;k)$.
\end{quote}

Restriction to the case $k\leq v/2$ entails no loss, since the map
$U\mapsto U^\perp$ (orthogonality being taken with respect to the
usual dot product) preserves the subspace distance and hence
identifies $q$-ary $(v,M,2\delta;k)$ subspace codes with
$(v,M,2\delta;v-k)$ subspace codes.

The Main Problem of Subspace Coding arose in connection with the
Koetter-Kschischang model of fault-tolerant network coding
\cite{koetter-kschischang08} (see \cite{kschischang12} for an
introduction), which
uses appropriate subspace codes to encode messages before
transmission over an ordinary network-coded (implemented by some form of
random linear network coding) packet network. Subspace codes with large
size and large minimum distance account for large transmission rate
and good error-correcting capabilities, respectively, and the
determination of the best such codes is thus of particular importance.

The Main Problem of Subspace Coding is akin to its classical
counterpart, the Main Problem of Algebraic Coding Theory, which asks
for the determination of the best linear codes over $\F_q$ relative to
the Hamming distance and forms the mathematical abstraction of the
engineering problem of finding the best point-to-point channel codes.
The Main Problem of Subspace Coding is much more difficult, however, 
since the set of all subspaces of $\F_q^v$ does not admit a group
structure compatible with the subspace metric. Hence there is no
suitable notion of ``linearity'' for subspace codes, which could be exploited.

In Finite Geometry language, a 
$q$-ary $(v,M,2\delta;k)$ subspace code is a set $\mathcal{C}$
of $M$ distinct ($k-1$)-flats in
the projective geometry $\PG(v-1,q)$ such that any ($k-\delta$)-flat
is contained in at most one member of $\mathcal{C}$ and some
($k-\delta-1$)-flat is contained in at least two 
members of $\mathcal{C}$. The Main Problem of Subspace Coding is
therefore equivalent to the packing problem for ($k-1$)-flats 
in $\PG(v-1,q)$ when these flats are identified with sets of ($t-1$)-flats
(all ($t-1$)-flats contained in the given ($k-1$)-flat), where
$t=k-\delta+1$. Thus it comes as no surprise that most of the known
results on the Main Problem, at the time of its advent,
had been obtained by Finite Geometers. With a few exceptions, 
these results pertain to the extremal case $\delta=k$ of pairwise
disjoint ($k-1$)-flats, so-called spreads or partial spreads, and
really satisfactory results were known only for the ``line'' case $k=2$. 

While some progress has been made since then, the Main Problem (unlike
its classical counterpart) is still widely open. Koetter and
Kschischang, in their seminal paper \cite{koetter-kschischang08},
used so-called maximum-rank-distance (MRD) codes, found earlier by
Delsarte \cite{delsarte78a}, Gabidulin \cite{gabidulin85} and Roth
\cite{roth91}, and a suitable lifting construction to
produce a good approximation to optimal subspace codes for general
parameter sets. Etzion and Silberstein \cite{etzion-silberstein09}
(cf.\ \cite{etzion-silberstein13} for the latest improvements)
introduced the echelon-Ferrers construction as a method to augment
lifted MRD (LMRD) codes by further subspaces, while keeping their minimum
distance. Further constructions of subspace codes for specific parameter
sets, improving on the general methods mentioned so far but usually
relying heavily on computer searches, can be found in
\cite{kohnert-kurz08,etzion-vardy11a,braun-reichelt14}. The numbers
$\smax_q(v,2\delta;k)$ have been determined exactly for $v\leq
5$ (all $q$) and in the binary case also for $v=6$. This includes
$\smax_q(4,4;2)=q^2+1$ (realized by a line spread in $\PG(3,q)$), 
$\smax_q(5,4;2)=q^3+1$ (the maximal size of a partial line spread in
$\PG(4,q)$; cf.\ \cite{beutelspacher75}), and $\smax_2(6,4;3)=77$
(realized by $5$ different isomorphism types of optimal subspace codes; cf.\ 
\cite{smt:fq11proc}).
Moreover, from the recent ground-breaking discovery of the first $2$-analogue
of a Steiner triple system in \cite{braun-etal13} it is also known
that $\smax_2(13,4;3)=1597245$.

Our contribution in this paper is a new approach to the construction of good
subspace codes, which has its origin in the observation made in
\cite{smt:fq11proc} that removing certain subcodes from an LMRD
code may result in the opportunity to add even more subspaces to
the expurgated LMRD code, thereby surpassing the size of the LMRD
code, as well as any other subspace code containing an LMRD code.

In contrast with the construction of a binary $(6,77,4;3)$ subspace
code of Type~A in~\cite{smt:fq11proc}, which used the smallest
possible set of removed codewords, we propose to remove a much larger
set of codewords from an LMRD code. As it turns out, our approach is
capable of constructing a binary $(7,329,4;3)$ subspace code,
equalizing the current record size and
providing an alternative construction of a code with $M=329$ for the
parameter set $q=2$, $v=7$, $k=3$, $d=4$; cf.\ \cite{braun-reichelt14}.

The main result is described in
Section~\ref{sec:newapp}. Section~\ref{sec:2-analogue} contains
information about a putative binary $(7,381,4;3)$ subspace code, whose
existence/non-existence is a famous unsolved problem, and the
subsequent two sections contain preparatory material for our main
result and related subspace code constructions.

Throughout the paper we will use basic concepts and terminology from
Finite Geometry. Readers are referred to
\cite{dembowski68,hirschfeld98,nova2011} for the relevant background
information and any unexplained terms.




\section{The Putative $2$-Analogue of the Fano 
  Plane}\label{sec:2-analogue} 

An easy double-counting argument yields the bound $M\leq 381$ for any
binary $(7,M,4;3)$ subspace code $\mathcal{C}$. If equality holds then
with $t=2$, $k=3$, $v=7$ every $t$-dimensional subspace of $\F_2^v$
must be contained in precisely one $k$-dimensional subspace (codeword)
of $\mathcal{C}$. Such a structure is known as a Steiner system over
$\F_2$, and in the particular case $(t,k,v)=(2,3,7)$ under
consideration has been named ``$2$-analogue of the \textsc{Fano} plane
$\PG(2,2)$'', since $\PG(2,2)$ is the unique ordinary Steiner system with these
parameters. Steiner systems over finite fields
(and, more generally, combinatorial designs
over finite fields) were introduced by Thomas \cite{thomas87}, who then
was able to show that a putative $2$-analogue of the Fano plane cannot
be constructed using the obvious idea of combining three plane orbits
(of size $127$) of a Singer group of $\PG(6,2)$.

In an attempt to construct such a $2$-analogue $\mathcal{C}$, one can
proceed as follows. Out of the $\gauss{7}{4}{2}=11811$ solids in $\PG(6,2)$
($4$-dimensional subspaces of $\F_2^7$), $15\times 381=5715$ should
contain a codeword of $\mathcal{C}$ and $6096$ should not contain a
codeword of $\mathcal{C}$. We now fix one such solid as the subspace
$S=(0,0,0,*,*,*,*)$ of $\F_2^7$ and count the number of codewords of
$\mathcal{C}$ meeting $S$ in a subspace of fixed dimension.

\begin{lemma}
  Suppose $\mathcal{C}$ is a putative $2$-analogue of the Fano plane
  and $a_i=\#\{U\in\mathcal{C};\dim(U\cap S)=i\}$ for $0\leq i\leq 3$
  is the so-called \emph{intersection vector of $S$ with respect to
    $\mathcal{C}$}. Then $(a_1,a_1,a_2,a_3)$ is either $(128,224,28,1)$ or
  $(136,210,35,0)$, depending on whether $S$ contains a codeword of
  $\mathcal{C}$ or not, respectively.
\end{lemma}
This lemma also follows from the general theory of intersection numbers
for subspace designs, as developed in \cite{kiermaier-pavcevic14}.
\begin{proof}
  We prove only the case where $S$ does not contain a codeword of
  $\mathcal{C}$, i.e.\ $a_3=0$. The other case is proved similarly.
  
  Each of the $\gauss{4}{2}{2}=35$ lines in $S$ is contained in
  exactly one codeword and the $35$ codewords obtained in this way are
  distinct, since the planes in $S$ spanned by two of the lines are
  not in $\mathcal{C}$. This gives $a_2=35$.
  Double-counting yields that each point of $\PG(6,2)$ is contained in
  $21$ codewords. Hence each of the $15$ points in $S$ must be
  contained in $7$ codewords meeting $S$ in a line and in $14$
  codewords meeting $S$ in a point. Consequently, $a_1=15\cdot 14=210$
  and $a_0=381-35-210=136$.
\end{proof}
Further information about a putative $2$-analogue $\mathcal{C}$ of the
Fano plane can be derived along these lines. For example, every
hyperplane ($5$-flat) of $\PG(6,2)$ must contain exactly $45$
codewords, every $4$-flat exactly $5$ codewords, and every
point-hyperplane flag $(p,H)$ must be incident with exactly $5$
codewords $U$ (i.e., $p\subset U\subset H$). All this is not
yet sufficient, however, to narrow down the number of possible
configurations, so that a computer search becomes feasible.

\section{Augmented LMRD codes}\label{sec:MRD}

Suppose $k,m,n$ are positive integers with $k\leq m\leq
n$.\footnote{The symbol '$k$' has a different meaning within
  this section.} An $(m,n,k)$ \emph{maximum rank distance (MRD)} code
over $\F_q$ is a set $\mathcal{A}=\{\mat{A}_1,...,\mat{A}_{q^{nk}}\}$
of $q^{nk}$ distinct matrices in $\F_q^{m\times n}$ having minimum
rank distance
$\rdist(\mathcal{A})=\min\bigl\{\rank(\mat{A}_i-\mat{A}_j);1\leq
i<j\leq q^{nk}\bigr\}=n-k+1$. An argument similar to that used in the
proof of the Singleton bound for ordinary block codes shows that
$q$-ary $(m,n,k=n-d+1)$ MRD codes in $\F_q^{m\times n}$ have maximum
size subject to the requirement $\rdist(\mathcal{A})\geq d$,
accounting for their name. According to
\cite{delsarte78a,gabidulin85,roth91}, MRD codes exist for all
admissible parameters $q,k,m,n$ and may be constructed using a
$q$-analogue of the familiar Reed-Solomon code construction (employing
linearized polynomials in place of ordinary polynomials).

As shown in \cite{koetter-kschischang08,silva-kschischang-koetter08},
the map $\lambda$
sending a matrix $\mat{A}\in\F_q^{m\times n}$ to the row space
of $(\imat_m|\mat{A})\in\F_q^{m\times(m+n)}$ satisfies
$\sdist\bigl(\lambda(\mat{A}),\lambda(\mat{B})\bigr)=2\rdist(\mat{A},\mat{B})$
(i.e. constitutes a ``scaled isometry'' with scale factor $2$). This
immediately gives that for any $q$-ary $(m,n,k)$ MRD code $\mathcal{A}$
the image $\lambda(\mathcal{A})$ forms a $q$-ary
$\bigl(m+n,q^{nk},2(n-k+1);m\bigr)$ subspace code, a so-called
\emph{lifted maximum rank distance (LMRD) code}.

We are interested in the case $q=2$, $(m,n,k)=(3,4,2)$, since an MRD
code with these parameters gives rise to a binary $(7,256,4;3)$ subspace
code, providing a good approximation to binary
optimal $(7,M,4;3)$ subspace codes. The standard MRD code with these
parameters is the ``Gabidulin code''

$\mathcal{G}=\{a_0x+a_1x^2;a_0,a_1\in\F_{16}\}$, viewed as a set of
$\F_2$-linear transformations $W\to\F_{16}$, $x\mapsto a_0x+a_1x^2$ on
a fixed $3$-dimensional $\F_2$-subspace $W$ of $\F_{16}$, which is
conveniently taken as the set of all elements $u\in\F_{16}$ of absolute
trace zero (i.e.\ $\trace(u)=\trace_{\F_{16}/\F_2}=u+u^2+u^4+u^8=0$). 
If $\F_{16}$ is constructed as $\F_2[\alpha]$ subject to
$\alpha^4+\alpha+1=0$, we have
$W=\{0,1,\alpha,\alpha^2,\alpha^4,\alpha^5,\alpha^8,\alpha^{10}\}$.
An explicit representation of $\mathcal{G}$ by binary $3\times 4$
matrices can then be obtained through fixing bases of $W$ and
$\F_{16}$ over $\F_2$ and using coordinates with respect to these
bases. 

In our implementation we have used $W=\langle
1,\alpha,\alpha^2\rangle$, $\F_{16}=\langle
1,\alpha,\alpha^2,\alpha^3\rangle$. Then the $4\times 4$ matrices
corresponding to $x\mapsto\alpha x$ and $x\mapsto x^2$ (which
determine all $256$ matrices of the $4\times 4$ matrix representation
of $\mathcal{G}$) are
\begin{equation*}
  \begin{pmatrix}
    0&1&0&0\\
    0&0&1&0\\
    0&0&0&1\\
    1&1&0&0
  \end{pmatrix}\quad\text{and}\quad
  \begin{pmatrix}
    1&0&0&0\\
    0&0&1&0\\
    1&1&0&0\\
    0&0&1&1
  \end{pmatrix},
\end{equation*}
respectively, and the $3\times 4$ matrix representation
of $\mathcal{G}$ is obtained by deleting the last row of all $256$
matrices. The rank of the $255$ nonzero $3\times 4$ matrices (and
hence the subspace distance between any two distinct $3\times 4$
matrices) is at least $2$, since a nonzero $\F_2$-linear transformation
$W\to\F_{16}$, $x\mapsto a_0x+a_1x^2$ has a kernel of dimension $\leq
1$. Applying $\lambda$ to the
$256$ selected $3\times 4$ matrices produces the required
$(7,256,4;3)$ LMRD code $\lambda(\mathcal{G})$
as an explicit set of $256$ generating matrices in
$\F_2^{3\times 7}$.
\begin{remark}
  \label{rmk:basis-free}
  Sometimes it is more convenient to work with the polynomials in
  $\mathcal{G}$ directly rather than with their representing
  matrices. A basis-independent representation of
  $\lambda(\mathcal{G})$ can be obtained as follows: Take
  $W\times\F_{16}\cong\F_2^7$ as the ambient vector space and the
  codewords of $\lambda(\mathcal{G})$ as the \emph{graphs} (in the sense of
  Real Analysis) of the linear maps induced by the polynomials in
  $\mathcal{G}$. In this representation the $256$ codewords of
  $\lambda(\mathcal{G})$ are
  \begin{equation*}
    G(a_0,a_1)=\bigl\{(x,a_0x+a_1x^2);x\in W\},\quad a_0,a_1\in\F_{16}.
  \end{equation*}
  Since $x\mapsto a_0x+a_1x^2$ is $\F_2$-linear, it is clear that each
  set $G(a_0,a_1)$ is a $3$-dimensional $\F_2$-subspace of
  $W\times\F_{16}$. Moreover, using coordinates in
  $W\times\F_{16}$ with respect to the ordered basis
  $(1,0)$, $(\alpha,0)$, $(\alpha^2,0)$, $(0,1)$, $(0,\alpha)$,
  $(0,\alpha^2)$, $(0,\alpha^3)$ identifies the spaces $G(a_0,a_1)$
  with the codewords of $\lambda(\mathcal{G})$ as introduced earlier.
\end{remark}


It has been shown in \cite{trautmann-rosenthal10} using the concept of
``pending dots'' that the lifted $(7,256,4;3)$ Gabidulin code
$\lambda(\mathcal{G})$ can be
augmented by $35$ further subspaces to a $(7,291,4;3)$ subspace code.
Any further enlargement is impossible, since the codewords of
$\lambda(\mathcal{G})$ cover all $7\times 256=1792$ lines disjoint
from the special solid $S=(0,0,0,*,*,*,*)$, so that additional
codewords must meet $S$ at least in a line, which clearly conflicts
with some of the $35$ codewords outside $\lambda(\mathcal{G})$ already
chosen.

We close this section by providing a different geometric construction
of the augmented $(7,291,4;3)$ subspace code, which seems to be worth
mentioning. It is known that $\PG(3,2)$ admits a line packing
(parallelism)
$\mathscr{S}=\{\mathcal{S}_0,\mathcal{S}_1,\mathcal{S}_2,\mathcal{S}_3,\mathcal{S}_4,\mathcal{S}_5,\mathcal{S}_6\}$,
i.e.\ the $35$ lines are partitioned into $7$ spreads $\mathcal{S}_i$,
each spread containing $5$ pairwise disjoint lines.\footnote{Such a
  line packing provides a solution to Kirkman's School Girl Problem;
  see e.g.\ \cite{mesner67}.}
Now take seven points $p_i$, $0\leq i\leq 6$, in $\PG(6,2)$ such that
the seven $4$-flats $\langle p_i,S\rangle$ are distinct (and hence
represent all the $4$-flats above $S$). Connect a line $L$ in $S$ to $p_i$
if $L\in\mathcal{S}_i$, and augment $\lambda(\mathcal{G})$ by the $35$
planes $\langle p_i,L\rangle$ obtained in this way. It is readily
checked that the resulting subspace code $\mathcal{C}$ of size $291$ still has
$\sdist(\mathcal{C})=4$.

Using the known ``cyclic'' description of $\mathscr{S}$ it is not hard
to produce a list of the $35$ additional codewords. In the
basis-independent representation one can take
$p_0=\F_2(1,0)=\F_2\times\{0\}$ and $\mathcal{S}_0$ to consist of
$\{0\}\times\alpha^j\F_4$ for $0\leq j\leq 4$, i.e.\ start with the
$5$ planes $E(0,j)=\F_2\times\alpha^j\F_4$, $0\leq j\leq 4$, and
generate the remaining codewords from these by applying the
linear permutation
$\sigma=(0)(\alpha^{14})(1,\alpha,\alpha^2,\alpha^4,\alpha^5,\alpha^{10},\alpha^8)
(\alpha^7,\alpha^{13},\alpha^9,\alpha^{12},\alpha^{11},\alpha^6,\alpha^3)$
of $\F_{16}$ simultaneously to both coordinates of $W\times\F_{16}$.
The resulting $35$ additional codewords are
\begin{equation*}
  E(i,j)=\sigma^i(\F_2)\times\sigma^i(\alpha^j\F_4),\quad 0\leq i\leq
  6,\;0\leq j\leq 4.
\end{equation*}
The choice of a point of the special form $p=\F_2(a,0)$ in every
additional plane is not mandatory, and in fact one could make $35$
arbitrary choices for the second coordinates of these points.

\section{Expurgating the LMRD Code First}\label{sec:expurgate}

From \cite{smt:fq11proc} we know that it makes sense to remove
certain sets of codewords from the $(7,256,4;3)$ code
$\lambda(\mathcal{G})$. The idea behind this 
approach is that the removal of $M_0$ codewords from
$\lambda(\mathcal{G})$ ``frees'' $7M_0$ lines disjoint from the
special solid $S$, which are no longer covered by the expurgated
subspace code, and hence can possibly be rearranged, four lines at a
time, into ``new planes'' $N$ of $\PG(6,2)$ meeting $S$ in a point. In
the best case, it will be possible to add $7M_0/4$ new planes to the
expurgated subspace code, resulting in a subspace code of size
$256-M_0+7M_0/4=256+3M_0/4$ that is superior to
$\lambda(\mathcal{G})$.

The results in \cite{smt:fq11proc} imply that all $7M_0$
free lines can be rearranged into new planes if the removed set of
codewords has
the form $\lambda(\mathcal{R})$ for some (disjoint) union
$\mathcal{R}=\biguplus_{s=1}^t(f_s+\mathcal{T})$ of cosets
of the following special $3$-dimensional $\F_2$-subspace $\mathcal{T}$
of $\mathcal{G}$:
\begin{equation*}
  \mathcal{T}=\{u^2x+ux^2;u\in W\}.
\end{equation*}
Any such choice of $\mathcal{R}$ uniquely determines $14t$ new planes meeting
$S$ in a point and covering, together with the codewords in the
expurgated subspace code, each line disjoint from $S$ exactly once. It
remains to be checked whether the new planes $N_i$ mutually satisfy
the subspace distance condition $\sdist(N_i,N_j)\geq 4$ if $i\neq
j$. Equivalently, if $N_i$ and $N_j$ pass through the same point $s\in S$
then $N_i\cap N_j=\{s\}$. (For new planes passing through distinct
points of $S$ there is no further restriction.)

Using the computer algebra package SAGE (\texttt{www.sagemath.org}),
we have checked how many cosets of $\mathcal{T}$ can be put into
$\mathcal{R}$ without violating the subspace distance condition for
the resulting new planes. It turned out that the maximum is $t=2$ and
$\mathcal{R}$ can be taken as $\{u^2x+ux^2;u\in\F_{16}\}$. This
results in a $(7,268,4;3)$ subspace code, which can be further
augmented by $35$ planes meeting $S$ in a line (using the method
described at the end of
Section~\ref{sec:MRD} with some specific choice of the points $p_i$)
to a $(7,303,4;3)$ subspace code.

A straightforward extension of the reasoning in \cite{smt:fq11proc}
shows that the $28$ new planes, obtained by rearranging the $112$ lines
in the planes $G(u^2,u)$, $u\in\F_{16}$, corresponding to $\mathcal{R}$,
meet $S$ in the points $\F_2(a^2b+ab^2)$ with $a,b\in W$ nonzero
and distinct, and have the basis-independent representation
\begin{equation*}
  N(Z,u)=\bigl\{(x,u^2x+ux^2+y);x\in Z,y\in\F_2(a^2b+ab^2)\bigr\},
\end{equation*}
where $Z=\langle a,b\rangle$ denotes one of the seven $2$-dimensional
$\F_2$-subspaces of $W$ und $u\in\F_{16}/Z$. 

From this we realized that the intersection points with $S$ of the new planes
determined by $\mathcal{R}$ are simply the 
points on $W$ (which forms a subplane of $S$), and thus
account for only $7$ of the $15$ points on
$S$. This is in contrast with the construction in \cite{smt:fq11proc},
which has $\dim(S)=3$ and new planes passing through every point of $S$.

Extending our scope to the ``rotated'' subspaces
$\mathcal{T}v=\bigl\{(u^2x+ux^2)v;u\in W,v\in\F_{16}^\times\}$, we
were able to overcome this restriction, but now had $15\times 32=480$
cosets to consider simultaneously. Moreover, cosets
$f+\mathcal{T}v$ for distinct $v$ need no longer be
disjoint, imposing an additional restriction.

The new problem can be viewed as a maximum clique problem in graph
theory.\footnote{A \emph{clique} in an undirected graph $G = (V, E)$
  is a subset $C\subseteq V$ of the vertex set such that any two
  vertices in $C$ are connected by an edge in $E$. A \emph{maximum
    clique} is a clique of the largest possible size $\#C$, called the
  \emph{clique number} of $G$.}  View
every coset $f+\mathcal{T}v$ as a vertex of an undirected graph
$G$. If two cosets are disjoint and their subspace lifts
have subspace distance at
least $4$ from each other, draw an edge between these two
vertices. The clique number of $G$ then gives 
the number of cosets $f+\mathcal{T}v$ we can put into
$\mathcal{R}$.

The clique number of $G$ turned out to be $4$ (again with
the help of SAGE), i.e. $4\times 8=32$ planes can be removed in
exchange for 56 new
planes, resulting in a $(7,280,4;3)$ subspace code. The
augmentation problem for this code can be modelled as a maximum clique
problem as well, and we found this time that $34$ further planes can
be added to produce a $(7,314,4;3)$ subspace code.


\section{The New Approach}\label{sec:newapp}

From Section~\ref{sec:2-analogue} we know that in order to qualify for
a putative $2$-analogue of the Fano plane, the ``removed set'' of
subspaces of $\lambda(\mathcal{G})$ should be much larger than those
considered in the previous section---about half the size of
$\mathcal{G}$. In the case where $S$ does not contain a codeword
(which can always be assumed by suitably changing the coordinate
system), we should remove $120=15\cdot 8$ planes (i.e.\ $15$ cosets
$f+\mathcal{T}v$)
from $\lambda(\mathcal{G})$ and replace these by $210=14\times 15$ new
planes ($14$ new planes through each point of $S$). A moment's
reflection shows that there is an obvious candidate for the
corresponding removed set of matrices
$\mathcal{R}$, viz.\ take
\begin{equation*}
  \mathcal{R}=\biguplus_{v\in\F_{16}^\times}(u^2x+ux^2+\mathcal{T})v,
  \quad\text{where
    $\trace(u)=1$}.
\end{equation*}
These $15$ cosets are disjoint, since $\mathcal{G}$ admits a
partition into $\{0\}$,
$\F_{16}^\times x$, $\F_{16}^\times x^2$ and the $15$ sets
$\bigl\{(u^2x+ux^2)v;u\in\F_{16}^\times\bigr\}$ with $v\in\F_{16}^\times$.

Using this set $\mathcal{R}$ as the removed set, it is at least 
conceivable that the resulting $7\times 120=4\times 210$ free lines
can be rearranged into $210$ new planes satisfying the subspace
distance condition. This condition is in fact quite easy to check,
since problems can occur only for the $14$ new planes passing through a
fixed point of $S$, and the $15$ sets of $14$ new planes determined in
this way are isomorphic as subspace codes (since multiplication by
$v\in\F_{16}^\times$, viewed as a Singer group acting on $S=(0,0,0|*,*,*,*)$,
identifies these sets).


The key question therefore is:
What is the size of the largest clique in one of these 14-sets of
new planes (say, the $14$ new planes through $p=\F_2(0,0,0|1,0,0,0)$)? 

Using again a maximum clique model, we found that
the clique number is $11<14$ (hence a putative $2$-analogue of the
Fano plane cannot be
constructed in this way), and the number of maximum cliques 
is $4$. This yields a new subspace code of size
$M=256-120+11\times15=301$, subject to further augmentation by
planes meeting $S$ in a line or being contained in $S$.
However, since there are $4$ choices for the $11$ new planes through
each point of $S$, the total number of new $(7,301,4;3)$ subspace
codes obtained in this way is $4^{15}=1073741824$. It is impossible to
check all these subspace codes for further augmentation in a reasonable
amount of time. Instead we used a randomized search method (checking several
thousands of cases) and found a maximum of $28$ planes that can be
added to some (in fact, many different) of the $4^{15}$ subspace codes,
resulting in a binary $(7,329,4;3)$ subspace code. This is our main
result and equalizes the record set in \cite{braun-reichelt14}.


\section{Conclusion}\label{conclusion}

We have outlined a new framework for the construction of good binary
$(v,M,4;3)$ subspace codes, which starts with a distinguished
($v-3$)-dimensional subspace $S$ of the ambient space $\F_2^v$ and
selects codewords based on their intersection dimension with
$S$. The subspace codes constructed do not contain lifted MRD codes
and hence are able to overcome the size limit imposed on such codes.
We have worked out the case $v=7$ in
detail and found that our framework is capable of yielding the largest
known subspace codes in this case.

Several challenging questions arise from our work. Is it possible to
construct the recently found $2$-analogue of a Steiner triple system,
a binary $(13,1597245,4;3)$ subspace code along these lines? For this
the distinguished subspace $S$ would be $10$-dimensional and, as one
can show, admit many different feasible intersection vectors
$(a_0,a_1,a_2,a_3)$. There is, however, a particular choice for $a_0$,
which is motivated by the example $v=7$ and determines the
intersection vector completely: $a_0=2^{19}+2^9=524\,800=2^9\times
5^2\times 41$, $a_1 = 916\,608=2^7\times 3\times 7\times 11\times 31$,
$a_2 = 152\,768=2^6\times 7\times 11\times 31$, $a_3 = 3069=3\times
11\times 31$. Thus $S$ would have to contain $3069=3\times 1023$
(three times the number of points in $S$) codewords, and
through each point $p\in S$ there would be $2^7\times 7=896$ new
planes meeting $S$ in $p$. Is such a construction of a 
binary $(13,1597245,4;3)$ subspace code actually possible?

Further questions are those for the largest $q$-ary $(7,M;4;3)$
subspace codes constructible by our method for prime powers $q>2$, and
for general lower bounds on the clique numbers of the graphs involved
that could be used to derive an infinite family of, say, binary
$(v,M,4;3)$ subspaces codes improving on the known general
constructions.

\section*{Acknowledgment}

The authors wish to thank Michael Kiermaier, University of Bayreuth,
for valuable discussions related to this paper and three reviewers for
their comments/suggestions on the initial submission.



\def\cprime{$'$}

\end{document}